\documentclass[12pt]{article}

\usepackage{latexsym}
\usepackage{amssymb}
\usepackage{amsfonts}
\usepackage{amsopn}

\newcommand{\ud}{\mathrm{d}}
\newcommand{\Id}{\mathrm{Id}}

\newtheorem{theorem}{Theorem}
\newtheorem{definition}{Definition}
\newtheorem{proposition}{Proposition}

\newtheorem{corollary}{Corollary}
\newtheorem{remark}{Remark}
\newtheorem{lemma}{Lemma}

\newcommand{\fn}{\footnote}

\def\bs{\vskip24pt}
\def\ms{\vskip12pt}

\def\nt{\noindent}

\def\title #1{\begin{center}
{\Large {\bf #1}}
\end{center}}
\def\author #1{\begin{center} {\large #1}
\end{center}}

\begin{document}

\begin{titlepage}
\def\thefootnote{\fnsymbol{footnote}}
\vspace*{1.1in}

\title{Individual Risk and Lebesgue Extension\\
\ms without Aggregate Uncertainty\fn{The authors are grateful to Bob Anderson, Jerry Keisler, Ali Khan, Rich McLean and Rabee Tourky
for the encouragement to us for working on the problem studied here. They also thank Peter Loeb for some useful discussions on the work of Fremlin.}}

\vskip .75em

\author{Yeneng Sun\footnote{Department of Economics,
National University of Singapore, 1 Arts Link, Singapore, 117570.\\
Email: ynsun@nus.edu.sg } and\, Yongchao Zhang\footnote{Department
of Mathematics, National University of Singapore, Science Drive 2,
Singapore, 117543.\\  Email: yongchao@nus.edu.sg}}

\vskip .75em

\centerline{First version: December 2007; this version: February
2008}

\vskip .75em \vskip 1.90em \baselineskip=.17in

\begin{abstract}
Many economic models include random shocks imposed on a large
number (continuum) of economic agents with individual risk. In this context, an exact law of
large numbers and its converse is presented in \cite{Sun-06} to
characterize the cancelation of individual risk via aggregation.
However, it is well known that the Lebesgue unit interval is not
suitable for modeling a continuum of agents in the particular
setting. The purpose of this note is to show that an extension of
the Lebesgue unit interval does work well as an agent space with
various desirable properties associated with individual risk.
\end{abstract}

\bigskip
\bs


\end{titlepage}

\setcounter{footnote}{0}

\section{Introduction}

Models with a continuum of agents are widely used in economics.
One often chooses to work with the the unit interval with the Lebesgue measure as the agent space.
However,  it was already noted by Aumann that the choice of the Lebesgue unit interval as a model for the agent space is of no particular significance and any atomless probability space is precisely what is needed to ensure that each individual agent has no influence.\fn{For this point, see p.~44 of \cite{Aumann}. For the discussion of various other formulations of negligible agents, see \cite{Khan-07}.}

Many economic models have also been based on a continuum of agents with
individual risk. Formally, a continuum of independent random
variables is used to model individual level random shocks imposed on
a large number of economic agents. The desirable result is an exact
law of large numbers which guarantees the cancelation of individual
risk at the aggregate level.\fn{See \cite{Sun-06} for many earlier
references on this. For some more recent applications of the law of
large numbers, see \cite{DGP-05}, \cite{DGP-07},
\cite{LR}, \cite{MP02}, \cite{MP05}, \cite{SY}, \cite{Weill}.} It is
shown in \cite{Sun-06} that a process measurable in a Fubini
extension is essentially pairwise independent if and only if it
satisfies the property of coalitional aggregate
certainty.\footnote{See Definitions \ref{df-epi} and
\ref{D:coalitional aggregate certainty} below respectively for the
precise meaning of essential pairwise independence and coalitional
aggregate certainty. The equivalence result is shown in Theorem 2.8
of \cite{Sun-06}.} The latter means that aggregation at the
coalitional level removes uncertainty.

Section 5 of \cite{Sun-06} considers the existence of a Fubini
extension that allows one to construct processes with essentially
pairwise independent random variables taking any given variety of
distributions. Many probability spaces can be used as the relevant
index space; for example, the big class of atomless Loeb probability
spaces (see \cite{LW}). One can also work with an index space based
on some atomless measure space on the unit interval $[0, 1]$.

However, unlike the case of a continuum of agents in a deterministic model,
it is well known that an economic model with an i.i.d. process based the classical
continuum product space and indexed by the Lebesgue unit interval
has the sample measurability problem.\fn{See \cite{Doob-37},
\cite{Doob-53}, \cite{Judd} and the detailed discussion in
\cite{Sun-06}.} Moreover, Corollary 4.3 of \cite{Sun-06} shows that
under the framework of a Fubini extension, almost all sample
functions of any process with essentially pairwise independent
random variables cannot be Lebesgue measurable. It is also pointed
out by Feldman and Gilles in Section 1 of \cite{Feldman and Gilles}
that a continuum of i.i.d. Bernoulli random variables indexed by the
Lebesgue unit interval cannot satisfy the the property of
coalitional aggregate certainty.

Since the Lebesgue unit interval is the simplest atomless probability space and it is
not suitable for modeling a continuum of agents with individual risk, a natural question is
whether one can find some extension of the Lebesgue unit interval as
the agent space with the desired property. The purpose of this note
is to provide a positive answer to that question. In particular, we
construct essentially pairwise independent processes measurable in a
Fubini extension, where the sample functions are measurable with
respect to some extension of the Lebesgue unit interval and the
random variables can take any given variety of distributions. It
follows immediately from the exact law of large numbers that the
type of result as mentioned in Section 1 of Feldman and Gilles
\cite{Feldman and Gilles} holds for some extension of the Lebesgue
unit interval. The point is that though the Lebesgue unit interval
fails to be an agent space modeling a continuum of agents with
individual risk, some extension of it does work.

The rest of the paper is organized as follows. Section \ref{Sec:
Basic results} presents some basic definitions and a previous
characterization result on the cancelation of individual risk via
aggregation. The main result is stated in section 3. The proofs are
given in the Appendix.

\section{Basics}\label{Sec: Basic results}

Let $(I, \mathcal{I}, \lambda)$ be an atomless probability space
which is used to model the agent space of many economic agents. In
our setting, it will be the parameter space for a process. Let
$(\Omega, {\cal F}, P)$ be a sample probability space, which models
the space of uncertain states of the world. A process $f$ from $I
\times \Omega$ to a complete separable metric space $X$ with Borel
$\sigma$-algebra~$\mathcal {B}$ is a mapping from $I \times \Omega$
to $X$ such that (1) for $\lambda$-almost all $i \in I$, the random
shock $f_i$ imposed on agent $i$ is a random variable defined on
$(\Omega, {\cal F}, P)$ whose distribution $P f_i ^{-1} $ on $X$ is
defined by $P f_i ^{-1} (B) = P [f_i ^{-1} (B)]$ for each $B \in \mathcal {B}$;
(2) for every $B \in \mathcal {B}$, the mapping
$i \mapsto P f_i ^{-1} (B)$ is $\mathcal {I}$-measurable.

The meaning of individual risk is that each individual agent is
allowed to have correlation with a negligible group of other agents.
This is formalized as the concept of  essential pairwise
independence.

\begin{definition} \label{df-epi} A process $f$ from $I \times \Omega$ to a
complete separable metric space $X$ is said to be essentially
pairwise independent if for $\lambda$-almost all $s \in I$, the
random variables $f_s$ and $f_i$ are independent for
$\lambda$-almost all $i \in I$.
\end{definition}

A desirable result for individual risk is its cancelation at the
aggregate level. That is, aggregation over a non-negligible group of
agents leads to no uncertainty. The following is a formal definition
of coalitional aggregate certainty.

\begin{definition} \label{D:coalitional aggregate certainty}
Let $f$ be a process from $I \times \Omega$ to a complete separable
metric space $X$. For any coalition $S$ (i.e., $S \in {\cal I}$ with
$\lambda(S) > 0$), let $f^S$ be the restriction of $f$ to $S \times
\Omega$, ${\cal I}^S = \{C \in {\cal I}: C \subseteq S\}$, and
$\lambda^S $ the probability measure rescaled from the restrictions
of $\lambda$ to ${\cal I}^S $. The process $f$ is said to satisfy
the property of coalitional aggregate certainty if for $P$-almost
all $\omega \in \Omega$, the sample function $f_\omega$ is ${\cal
I}$-measurable, and for each coalition $S$, the empirical (sample)
distribution $\lambda (f^S_\omega)^{-1} (\cdot)$ is $\int_S
P f_i ^{-1} (\cdot) d\lambda^S$ for $P$-almost all $\omega \in \Omega$.
\end{definition}

When $X$ is the real line $\mathbb{R}$ and the random variables
$f_i$ are i.i.d. with a common distribution function $F$,
coalitional aggregate certainty means that for each coalition $S$,
the empirical distribution function $F^S_\omega$ generated by the
restricted sample function $f^S_\omega$ is $F$ for $P$-almost all
$\omega \in \Omega$. It is easy to construct examples of a continuum
of independent random variables with this aggregation property for
the grand coalition $I$ or for all the coalitions; see the
discussion in Section 6.3 of \cite{Sun-06}, and \cite{Anderson},
\cite{Green}, \cite{Judd}. However, one can also construct other
examples of a continuum of independent random variables whose sample
functions may not be measurable, or behave in a very ``strange" way.
In fact, for an i.i.d. process based on the usual continuum product
via the Kolmogorov construction, one can obtain the absurd claim
that almost all sample functions are essentially equal to an
arbitrarily given function $h$ on the index space (see Proposition
6.1 of \cite{Sun-06}); and thus, the sample distribution can be
undefined or completely arbitrary (see also \cite{Doob-37} and
\cite{Judd}).

The main difficulty for working with an essentially pairwise
independent process $f$ is that if it is jointly measurable with
respect to the usual product $\sigma$-algebra ${\cal I} \otimes
{\cal F}$, then the random variables $f_i$ are essentially constant
for almost all $i \in I$ (see Proposition 2.1 of \cite{Sun-06}).
Consequently, the usual product probability space $(I \times \Omega,
{\cal I} \otimes {\cal F}, \lambda \otimes P)$ will (typically) be
inadequate to prove any meaningful result on no aggregate
uncertainty. As shown in \cite{Sun-06}, a simple way to resolve this
problem is to work with an extension of the usual product
probability space that retains the Fubini property.

\begin{definition}\label{D:Fubini extension}
Let $(I \times \Omega, {\cal I} \otimes {\cal F}, \lambda \otimes
P)$ be the usual product probability space of the two probability
spaces $(I, \mathcal{I}, \lambda)$ and $(\Omega, {\cal F}, P)$. A
probability space $(I \times \Omega,
 {\cal W}, Q)$ extending $(I \times \Omega, {\cal I}
 \otimes {\cal F}, \lambda \otimes P)$ is said to be a \emph {Fubini
extension} if for any real-valued $Q$-integrable function $f$ on
$(I\times \Omega, {\cal W})$,

(1) the two functions $f_i$ and $f_\omega$ are integrable
respectively on $(\Omega, {\cal F}, P)$ for $\lambda$-almost all $i
\in I$, and on $(I, {\cal I}, \lambda)$ for $P$-almost all $\omega
\in \Omega$;

(2) $\int_\Omega f_i dP$ and $\int_I f_\omega dP$ are integrable
respectively on $(I, {\cal I}, \lambda)$ and $(\Omega, {\cal F},
P)$, with $ \int_{I \times \Omega} f dQ = \int_I \left( \int_\Omega
f_i dP \right) d\lambda = \int_\Omega \left( \int_I f_\omega
d\lambda \right) dP$.\fn{The classical Fubini Theorem is only stated
for the usual product measure spaces. It does not apply to
integrable functions on $(I \times \Omega, {\cal W}, Q)$ since these
functions may not be ${\cal I} \otimes {\cal F}$-measurable.
However, the conclusions of that theorem do hold for processes on
the enriched product space $(I \times \Omega, {\cal W}, Q)$ that
extends the usual product.}

To reflect the fact that the probability space $(I \times \Omega,
{\cal W}, Q)$
 has $(I, {\cal I}, \lambda)$ and $(\Omega, {\cal F}, P)$ as its marginal
spaces, as required by the Fubini property, it will be denoted by
$(I
 \times \Omega, {\cal I} \boxtimes {\cal F}, \lambda \boxtimes P)$.
\end{definition}

The following result is shown in Theorem 2.8 of \cite{Sun-06}. It
indicates that the framework of Fubini extension does deliver the
desired exact law of large numbers which guarantees the cancelation
of individual risk at the aggregate level.

\begin{lemma}\label{L:equivalence of fubini extension}
Let $f$ be a measurable process from a Fubini extension $(I\times \Omega, \mathcal{I}\boxtimes
\mathcal{F}, \lambda\boxtimes P)$ to a complete separable metric space $X$. Then $f$ satisfies the property of coalitional aggregate certainty
if and only if $f$ is essentially pairwise independent.
\end{lemma}

\section{The main result}\label{sec:Main
result}

Let $L=[0,1]$, $\mathcal{L}$ the $\sigma$-algebra of Lebesgue
measurable sets, and $\eta$ the Lebesgue measure defined on
$\mathcal{L}$. The Lebesgue unit interval is simply
$(L,\mathcal{L},\eta)$.

We are now ready to state the main result of this paper, which shows
that some extension $(I,\mathcal{I},\lambda)$ of the Lebesgue unit
interval $(L,\mathcal{L},\eta)$ can be used as the agent space
modeling a continuum of agents with individual risk in a very
general setting.\fn{$(I,\mathcal{I},\lambda)$ extending
$(L,\mathcal{L},\eta)$ means that $I=L =[0,1]$, $\mathcal{I}$
contains $\sigma$-algebra $\mathcal{L}$ of Lebesgue measurable sets,
and $\lambda$ extends the Lebesgue measure $\eta$ on $\mathcal{L}$.}
In particular, we show the existence of essentially pairwise
independent processes measurable in a Fubini extension, where the
sample functions are measurable with respect to the extended
Lebesgue interval $(I,\mathcal{I},\lambda)$ and the random variables
can take any given variety of distributions.

\begin{theorem} \label{T: main theorem} Let $I$ be the unit interval $[0,1]$ and
$X$ be a complete separable metric space. There exists a probability
space $(I,\mathcal{I},\lambda)$ extending the Lebesgue unit
interval, a probability space $(\Omega,\mathcal{F},P)$, and a Fubini
extension
$(I\times\Omega,\mathcal{I}\boxtimes\mathcal{F},\lambda\boxtimes P)$
such that for any measurable mapping $\varphi$ from $(I, {\cal I},
\lambda)$ to the space ${\cal M}(X)$ of Borel probability
measures\fn{${\cal M}(X)$ is endowed with the topology of weak
convergence of measures. The measurability of $\varphi$ is equivalent to the measurability of the mappings
$i \mapsto \varphi(i) (B) $ for $B \in \mathcal {B}$.} on $X$, there is a ${\cal I} \boxtimes
{\cal F}$-measurable process $f$ from $I \times \Omega$ to $X$ such
that the random variables $f_i$ are essentially pairwise
independent, and the distribution $P f_i^{-1}$ is the given
distribution $\varphi(i)$ for $\lambda$-almost all $i \in I$.
\end{theorem}

The following corollary on the special case of an i.i.d. process is obvious.

\begin{corollary} \label{cor-main}
Let $I$ and $X$ be as in Theorem \ref{T: main theorem}. There exists
a probability space $(I,\mathcal{I},\lambda)$ extending the Lebesgue
unit interval, a probability space $(\Omega,\mathcal{F},P)$, and a
Fubini extension
$(I\times\Omega,\mathcal{I}\boxtimes\mathcal{F},\lambda\boxtimes P)$
such that for any Borel probability measure $\tau$ on $X$, there is
a ${\cal I} \boxtimes {\cal F}$-measurable process $f$ from $I
\times \Omega$ to $X$ such that the random variables $f_i$ are
essentially pairwise independent with common distribution $\tau$.
\end{corollary}

\begin{remark} By Lemma \ref{L:equivalence of fubini extension}, essential
pairwise independence implies coalitional aggregate certainty. Thus
the processes in Theorem \ref{T: main theorem} and  Corollary
\ref{cor-main} satisfy coalitional aggregate certainty. For the special case that $X$ is the real line $\mathbb{R}$ and the random variables
$f_i$ are i.i.d. with a common distribution function $F$ as in Corollary \ref{cor-main}, we obtain that for each coalition $S$,
the empirical distribution function $F^S_\omega$ generated by the
restricted sample function $f^S_\omega$ is $F$ for $P$-almost all
$\omega \in \Omega$. Thus, the type of result as mentioned in Section 1 of Feldman and
Gilles  \cite{Feldman and Gilles} holds for some extension of the
Lebesgue unit interval.
\end{remark}

\section{Appendix}\label{sec: appendix}

In this appendix, the unit interval $[0,1]$ will have a different
notation in a different context. Recall that $(L,\mathcal{L},\eta)$
is the Lebesgue unit interval. We shall often work with the case
that the target space $X$ is the unit interval $[0, 1]$ with uniform
distribution $\mu$. Here $\mu$ is simply the Lebesgue measure
defined on Borel $\sigma$-algebra ${\cal B}$ of $[0,1]$. Note that
the Lebesgue measure defined on the Lebesgue $\sigma$-algebra
$\mathcal{L}$ is denoted by $\eta$.

The following result is Proposition 5.6 of \cite{Sun-06}.

\begin{lemma}\label{L:egz of with fubini extension}
There is an atomless probability space $(K,\mathcal{K},\kappa)$ with
$K=[0,1]$, a probability space $(\Omega, \mathcal{F},P)$, a Fubini
extension $(K\times \Omega,\mathcal{K}\boxtimes
\mathcal{F},\kappa\boxtimes P)$, and a $\mathcal{K}\boxtimes
\mathcal{F}$-measurable process $g$ from $K \times \Omega$ to
$[0,1]$ such that the random variables $g_k$ are pairwise
independent and identically distributed (i.i.d.) with common uniform
distribution $\mu$ on $[0, 1]$.
\end{lemma}

In addition, in Proposition 5.6 of \cite{Sun-06}, the sample
probability space $(\Omega, \mathcal{F},P)$ is an extension of the
usual continuum product; the index space $(K,\mathcal{K},\kappa)$ is
obtained from a Loeb probability space via a bijection. As shown in
Corollary 3 of \cite{KS}, $(K,\mathcal{K},\kappa)$ is not an
extension of the Lebesgue unit interval $(L,\mathcal{L},\eta)$.
However, as mentioned earlier, the purpose of this paper is to
obtain some extension of the Lebesgue unit interval as an agent
space with various desirable properties associated with individual
risk.

The following is essentially taken from Lemma 419I of
Fremlin~\cite{Fremlin 4}. Its proof is based on the Transfinite
Induction.

\begin{lemma}\label{L:Fremlin 419I}
There is a disjoint family ${\cal C} = \{C_k: k \in K = [0, 1]\}$ of
subsets of $L = [0,1]$ such that $\bigcup_{k \in K} C_k = L$, and
for each $k \in K$, $\eta_*(C_k)=0$ and $\eta^*(C_k)=1$, where
$\eta_*$ and $\eta^*$ are the respective inner and outer measures of
the Lebesgue measure $\eta$.
\end{lemma}

The original version of Lemma 419I of Fremlin~\cite{Fremlin 4} does
not require $\bigcup_{k \in K} C_k = L$. Suppose $\bigcup_{k \in K}
C_k \ne L$; let $B = L \setminus \bigcup_{k \in K} C_k$. Since the
cardinality of $B$ is at most the cardinality of $[0, 1]$, we can
redistribute at most one point of $B$ into each $C_k$ in the family
${\cal C}$.

As in the proof of Lemma 521P(b) of~\cite{Fremlin 5}, we define a
subset $C$ of $L\times K$ by letting $C = \{(l,k) \in L\times K:
l\in C_k, k\in K\}$. Let $(L\times K, \mathcal{L}\otimes\mathcal{K},
\eta\otimes\kappa)$ be the usual product probability space. For any
$\mathcal{L}\otimes\mathcal{K}$-measurable set $U$ that contains
$C$, $C_k \subseteq U_k$ for each $k \in K$. The Fubini property of
$\eta\otimes\kappa$ implies that for $\kappa$-almost all $k \in K$,
$U_k$ is $\mathcal{L}$-measurable, which means that $\eta(U_k) =1$
(since $\eta^*(C_k)=1$). Since $\eta\otimes\kappa (U) = \int_K
\eta(U_k) \ud \kappa$, we have $\eta\otimes\kappa (U) = 1$.
Therefore, the $\eta\otimes\kappa$-outer measure of $C$ is one.

Since the $\eta\otimes \kappa$-outer measure of $C$
is one, the method in \cite{Doob-53} (see p. 69) can be used to extend $\eta \otimes \kappa$ to a measure $\gamma$ on the $\sigma$-algebra $\mathcal{U}$
generated by the set $C$ and the sets in ${\cal L \otimes K}$ with $\gamma (C) = 1$. It is easy to see that $\mathcal{U} = \left\{ (U_1 \cap  C)\cup(U_2 \backslash C): U_1,U_2\,\in \,\mathcal{L} \otimes
\mathcal{K} \right\}$, and $\gamma[(U_1\cap C)\cup(U_2 \backslash C)]
 =\eta\otimes\kappa(U_1)$ for any measurable sets $U_1,U_2\,\in
\,\mathcal{L}\otimes\mathcal{K}$. Let $\mathcal{T}$ be the $\sigma$-algebra $\{U\cap C:
U\in\mathcal{L}\otimes\mathcal{K} \}$, which is the collection of all the measurable subsets of $C$ in $\mathcal{U}$. The restriction of $\gamma$ to $(C, \mathcal{T})$ is still denoted by $\gamma$. Then, $\gamma(U\cap C) = \eta\otimes\kappa\,(U)$, for every measurable set
$U\in\mathcal{L}\otimes\mathcal{K}$. Note that $(L \times K, {\cal U}, \gamma)$ is an extension of $(L\times K, {\cal L \otimes K}, \eta
\otimes \kappa)$.

Consider the projection mapping $\pi:L\times K\rightarrow L$ with
$\pi(l,k) = l$. Let $\psi$ be the
restriction of $\pi$ to $C$. Since the family ${\cal C}$ is a
partition of $L=[0,1]$, $\psi$ is a bijection between $C$ and $L$.
It is obvious that $\pi$ is a measure-preserving mapping from
$(L\times K, \mathcal{L}\otimes\mathcal{K}, \eta\otimes\kappa)$ to
$(L,\mathcal{L},\eta)$, i.e., for any $B \in \mathcal{L}$, $\pi^{-1}
(B) \in \mathcal{L}\otimes\mathcal{K}$ and
$\eta\otimes\kappa[\pi^{-1} (B)] = \eta(B)$; and thus $\pi$ is a measure-preserving mapping from $(L \times K, {\cal U}, \gamma)$ to $(L,\mathcal{L},\eta)$.
Since $\gamma(C) = 1$, $\psi$ is a
measure-preserving mapping from $(C, \mathcal{T}, \gamma)$ to
$(L,\mathcal{L},\eta)$, i.e., $\gamma[\psi^{-1}(B)] = \eta(B)$ for
any $B \in \mathcal{L}$.

To introduce one more measure structure on $[0, 1]$, we shall also
denote it by $I$.  Let $\mathcal{I}$ be the $\sigma$-algebra $\{S
\subseteq I: \psi^{-1} (S) \in \mathcal{T}\}$.
Define a set function $\lambda$ on $\mathcal{I}$ by letting
$\lambda(S)= \gamma[\psi^{-1}(S)]$ for each $S \in {\cal I}$. Since $\psi$ is a bijection,
$\lambda$ is a well-defined probability measure on $(I,
\mathcal{I})$. Moreover, $\psi$ is also an isomorphism from $(C,
\mathcal{T}, \gamma)$ to $(I, \mathcal{I}, \lambda)$. Since $\psi$
is a measure-preserving mapping from $(C, \mathcal{T}, \gamma)$ to
$(L,\mathcal{L},\eta)$, it is obvious that $(I, \mathcal{I},
\lambda)$ is an extension of the Lebesgue unit interval
$(L,\mathcal{L},\eta)$.\fn{A similar observation was made in $\S6$
of~\cite{Podczeck 08} for the case of the product of the Lebesgue unit interval and the space $\{0, 1\}^\alpha$ with the cardinality $\alpha$ between the cardinality $c$ of the continuum and $2^c$. Thus the cardinality of $\{0, 1\}^\alpha$ is at least $2^c$ while the cardinality of our space $K$ in Lemma \ref{L:egz of with fubini extension} is $c$.}

As we have seen, based on the constructions as used in Lemma 419I of
Fremlin~\cite{Fremlin 4} and Lemma 521P(b) of~\cite{Fremlin 5}, it
is rather straightforward to construct $(I, \mathcal{I}, \lambda)$.
To prove Theorem \ref{T: main theorem}, the key part is to construct
a Fubini extension and essentially pairwise independent measurable
processes whose random variables take any variety of distributions.
We shall first consider a special case of Theorem \ref{T: main
theorem} below.

\begin{proposition}\label{P:rich product space.}
There is a Fubini extension
$(I\times\Omega,\mathcal{I}\boxtimes\mathcal{F},\lambda\boxtimes P)$
and an essentially pairwise independent process $f:I\times
\Omega\rightarrow [0,1]$ such that $f$ is
$\mathcal{I}\boxtimes\mathcal{F}$-measurable, and for each $i \in
I$, the distribution of the random variable $f_i$ is the uniform
distribution $\mu$ on $[0, 1]$.
\end{proposition}

\nt \textbf{Proof:}  We construct the process
$f:I\times\Omega\rightarrow [0,1]$ in three steps.\\
\nt $\textbf{Step 1.}$ Based on the process $g$ and the Fubini
extension $(K\times \Omega,\mathcal{K}\boxtimes
\mathcal{F},\kappa\boxtimes P)$ in Lemma~\ref{L:egz of with fubini
extension}, we construct a new process $G$ from the triple product
space $L\times K\times\Omega$ to $[0, 1]$ with $G(l,k,\omega) =
g(k,\omega)$ for each $(l, k, \omega) \in L\times K\times\Omega$. Here the index space is augmented to the usual product
space $(L\times K,\mathcal{L}\otimes\mathcal{K},\eta\otimes\kappa)$
while the sample space remains $(\Omega, \mathcal{F},P)$.

For each $(l,k) \in L \times K$, $G_{(l,k)} = g_k$ is a random
variable on the sample space with common uniform distribution $\mu$
on $[0,1]$. Moreover, the process $G$ is essentially pairwise
independent. In fact, for any $(l_0,k_0)\in L\times K$, if $k\neq
k_0 $,  $g_{k_0}$ and $g_k$ are independent random variables, so are
the random variables $G_{(l_0,k_0)} = g_{k_0}$ and $G_{(l,k)} =
g_k$. It is obvious that the subset $\{(l,k)\in L\times K: k\neq
k_0\}$ has full $\eta\otimes\kappa$-measure.

Now consider the usual product space $(L\times K\times\Omega,
\mathcal{L}\otimes(\mathcal{K}\boxtimes\mathcal{F}),\eta\otimes(\kappa\boxtimes
P))$ of the Lebesgue unit interval $(L,\mathcal{L},\eta)$ with the
Fubini extension $(K\times \Omega,\mathcal{K}\boxtimes
\mathcal{F},\kappa\boxtimes P)$. Note that the process $G$ is
$\mathcal{L}\otimes(\mathcal{K}\boxtimes\mathcal{F})$-measurable
because $g$ is $\mathcal{K}\boxtimes\mathcal{F}$-measurable. Next we
claim that it is a Fubini extension of the usual triple product
space $((L\times K)\times\Omega,
(\mathcal{L}\otimes\mathcal{K})\otimes\mathcal{F}),
(\eta\otimes\kappa)\otimes P)$.

To show the Fubini property on the extended space, we adapt a proof
analogous to the usual Fubini Theorem.\fn{See, for example, p.~308
of~\cite{Royden}. Similar adaption of the idea has been used in
\cite{HS-06} to prove the one-way Fubini property.} Let $V\subseteq
L^1(\eta\otimes(\kappa\boxtimes P))$ be the set of all
$\eta\otimes(\kappa\boxtimes P)$-integrable function $h$ satisfying
the Fubini property. That is, (1) $h_{(l,k)}$ is integrable on
$(\Omega, \mathcal{F},P)$ for $\eta\otimes\kappa$-almost all
$(l,k)\in L\times K$ and $h_{\omega}$ is integrable on $(L \times
K,\mathcal{L}\otimes \mathcal{K}, \eta\otimes\kappa )$ for
$P$-almost all $\omega\in\Omega$; (2) $\int_{\Omega}h_{(l,k)}\ud P$
and $\int_{L\times K} h_{\omega}\ud \eta\otimes\kappa$ are
integrable respectively on $(L\times
K,\mathcal{L}\otimes\mathcal{K},\eta\otimes\kappa )$ and
$(\Omega,\mathcal{F},P)$; (3) $\int_{\tiny{L\times K \times \Omega}}
h\, \ud \,(\eta\otimes(\kappa\boxtimes P)) = \int_{\tiny{L\times
K}}(\int_{\Omega}h_{(l,k)}\,\ud P)\ud \eta\otimes\kappa
=\int_{\Omega}(\int_{\tiny{L\times K}}
h_{\omega}\,\ud\eta\otimes\kappa)\ud P.$

We shall first show that the set $V$ contains all the indicator
functions of the measurable sets in
$\mathcal{L}\otimes(\mathcal{K}\boxtimes\mathcal{F})$. Let
$\mathcal{D}$ be the collection of all
$\mathcal{L}\otimes(\mathcal{K}\boxtimes\mathcal{F})$-measurable
sets $D$ such that its indicator function $1_D$ (which takes value
$1$ in $D$ and $0$ outside) is in $V$.

Consider $D$ to be a measurable rectangle $B\times W$ for $B \in
{\cal L}$ and $ W\in \mathcal{K\boxtimes F}$. The section $D_\omega$
is $B\times W_{\omega}$. By the Fubini property associated with
$\kappa\boxtimes P$, $B\times W_{\omega}$ is in
$\mathcal{L}\otimes\mathcal{K}$ for $P$-almost all $\omega\in
\Omega$. The measure of $D_\omega$ is $\eta\otimes\kappa (D_\omega)
= \eta(B) \kappa(W_{\omega})$, which is $P$-integrable with integral
$\eta(B) \int_\Omega \kappa(W_{\omega})\, \ud P$. Similarly, the
section $D_{(l,k)}$ is $W_k$ if $l \in B$, and empty if $l \notin
B$, which is in ${\cal F}$ for $\kappa$-almost all $k\in K$. The
measure of $D_{(l,k)}$ is $1_B(l) P(W_k)$, which is
$\eta\otimes\kappa$-integrable with integral $\eta(B) \int_K
P(W_k)\ud \kappa$. By the Fubini property associated with
$\kappa\boxtimes P$ again,
$$\int_\Omega \eta(B)
\kappa(W_{\omega}) \, \ud P =\eta(B) \int_\Omega  \kappa(W_{\omega})
\, \ud P = \eta(B) \int_K P(W_k)\ud \kappa = \eta(B)
(\kappa\boxtimes P) (W),$$ which means that $$\int_\Omega
\eta\otimes\kappa (D_\omega)\, \ud P= \int_{L \times K}
P(D_{(l,k)})\, \ud \eta\otimes\kappa = (\eta\otimes(\kappa\boxtimes
P)) (D).$$ Hence, $B\times W \in \mathcal{D}$.

Next, we show that the collection $\mathcal{D}$ is also a Dynkin (or
$\lambda$-) system on $L\times K\times \Omega$. Indeed, it is
obvious that (i) $L\times K\times \Omega\in \mathcal{D}$; (ii) if
$D,D'\in \mathcal{D}$ and $D'\subseteq D$, then $D-D'\in\mathcal{D}$
because $1_{D-D'}=1_D-1_{D'}$ and Fubini property is closed under
linear combination; (iii) if $D^n$ is an increasing sequence of sets
in $\mathcal{D}$, then $1_{D^n}$ is an increasing sequence of
functions with limit $1_{\cup^{\infty}_{n=1}D^n}$, thus
$\cup^{\infty}_{n=1}D^n\in\mathcal{D}$ according to the Monotone
Convergence Theorem (see the proof for a general sequence of
functions below). Since the collection of measurable rectangles of
the form $B\times W$ for $B \in {\cal L}$ and $ W\in
\mathcal{K\boxtimes F}$ is $\pi$-system (i.e., closed under finite
intersections) and generates
$\mathcal{L}\otimes(\mathcal{K}\boxtimes\mathcal{F})$, Dynkin's
$\pi$-$\lambda$ Theorem (see p.~277 of \cite{cA} or p.~24
of~\cite{Durrett}) implies that
$\mathcal{D}=\mathcal{L}\otimes(\mathcal{K}\boxtimes\mathcal{F})$.
Hence $V$ contains all the indicator functions of the measurable
sets in $\mathcal{L}\otimes(\mathcal{K}\boxtimes\mathcal{F})$.

As mentioned above, the set $V$ is closed under linear combinations.
In particular, $V$ contains all the measurable simple functions and
the difference between any two members. Note that each
$\eta\otimes(\kappa\boxtimes P)$-integrable function is the
difference between two non-negative integrable functions and each
non-negative integrable function is the pointwise limit of an
increasing sequence of non-negative simple functions. So we only
need to show that for any increasing sequence of non-negative
functions in $V$ with an integrable pointwise limit, the limit
function also belongs to $V$.

Now let $h\in L^1(\eta\otimes(\kappa\boxtimes P))$, and
$\{h^n\}_{n=1}^\infty$ be an increasing sequence of non-negative
functions in $V$ with pointwise limit $h$ (to be denoted by $h^n
\uparrow h$). By the Monotone Convergence Theorem (see
\cite{Royden}), $$\lim_{n\to \infty} \int_{\tiny{L\times K \times
\Omega}} h^n\, \ud \,(\eta\otimes(\kappa\boxtimes P)) =
\int_{L\times K \times \Omega} h\, \ud
\,(\eta\otimes(\kappa\boxtimes P)).$$ Since $h^n$ satisfies the
Fubini property, we know that $h^n_\omega$ is
$\eta\otimes\kappa$-integrable for $P$-almost all $\omega\in\Omega$.
For each $\omega \in \Omega$, $h^n_\omega \uparrow h_\omega$. The
Monotone Convergence Theorem implies that for $P$-almost all
$\omega\in\Omega$, $\int_{L\times K}h^n_{\omega}\, \ud
\eta\otimes\kappa \uparrow \int_{L\times K}h_{\omega}\, \ud
\eta\otimes\kappa$. By the Fubini property of $h^n$ again,
$\int_{K\times L}h^n_{\omega}\, \ud \eta\otimes\kappa$ is
$P$-integrable. Hence, the Monotone Convergence Theorem can be
applied again to obtain that
$$\lim_{n\to \infty} \int_{\Omega}\left(\int_{L\times K}
h^n_{\omega}\,\ud\eta\otimes\kappa\right)\ud P =
\int_{\Omega}\left(\int_{L\times K}
h_{\omega}\,\ud\eta\otimes\kappa\right)\ud P.$$ Since
$\int_{\Omega}\left(\int_{L\times K}
h^n_{\omega}\,\ud\eta\otimes\kappa\right)\ud P =  \int_{L\times K
\times \Omega} h^n\, \ud \,(\eta\otimes(\kappa\boxtimes P))$, we
have $$\int_{\Omega}\left(\int_{L\times K}
h_{\omega}\,\ud\eta\otimes\kappa\right)\ud P =  \int_{L\times K
\times \Omega} h\, \ud \,(\eta\otimes(\kappa\boxtimes P)).$$ The
other half of the Fubini property for $h$ can be proved in a similar
way. Hence $h \in V$.

Therefore, we show that $V= L^1(\eta\otimes(\kappa\boxtimes P))$,
which means that the extended space $(L\times K\times
\Omega,\mathcal{L}\otimes(\mathcal{K}\boxtimes\mathcal{F}),\eta\otimes(\kappa\boxtimes
P))$ is a Fubini extension.

\medskip

\nt $\textbf{Step 2.}$ Now consider a new process $F$ from $C\times
\Omega$ to $[0,1]$ with $F$ being the restriction $G|_{C\times
\Omega}$ of $G$ to $C\times \Omega$, where $G$ is the process in
Step 1. The index probability space is restricted to
$(C,\mathcal{T},\gamma)$ from the the product space $(L\times K,
\mathcal{L}\otimes\mathcal{K}, \eta\otimes\kappa)$, and the sample
space $(\Omega, {\cal F}, P)$ remains the same as in Step 1.

It is clear that for any $(l, k) \in C$, $F_{(l, k)}$ is a random
variable on the sample space with uniform distribution $\mu$ on
$[0,1]$. Moreover, $F$ is an essentially pairwise independent
process. In fact, for any $(l_0,k_0)\in C\subseteq L\times K$, the
random variables $F_{(l,k)}=G_{(l,k)}$ and
$F_{(l_0,k_0)}=G_{(l_0,k_0)}$ are independent for any $(l, k) \in
\{(l,k)\in L\times K : k\neq k_0\} \cap C $, note that
$$\gamma[ \{(l,k)\in L\times K: k\neq k_0\} \cap C]=\eta\otimes
\kappa [\{(l,k)\in L\times K: k\neq k_0\}]=1.$$

Recall that $(L\times K\times
\Omega,\mathcal{L}\otimes(\mathcal{K}\boxtimes\mathcal{F}),\eta\otimes(\kappa\boxtimes
P))$ is shown to be a Fubini extension in Step 1. We shall prove
that $C \times \Omega$ has $\eta\otimes(\kappa\boxtimes P)$-outer
measure one. Let $D$ be any measurable set in $\mathcal{L} \otimes
(\mathcal{K} \boxtimes \mathcal{F})$ that contains $C \times
\Omega$. Then, for each $\omega \in \Omega$, $C \subseteq D_\omega$.
By the Fubini property associated with $\eta\otimes(\kappa\boxtimes
P)$, we have for $P$-almost all $\omega\in\Omega$, $D_\omega \in
\mathcal{L}\otimes\mathcal{K}$, and hence $\eta\otimes \kappa
(D_\omega)=1$ since $C$ has $\eta\otimes\kappa$-outer measure one.
By the Fubini property associated with $\eta\otimes(\kappa\boxtimes
P)$ again, $\eta\otimes(\kappa\boxtimes P)(D)=\int_\Omega
\eta\otimes \kappa (D_\omega)\,\ud P = 1$.

Based on the Fubini extension $(L\times K\times
\Omega,\mathcal{L}\otimes(\mathcal{K}\boxtimes\mathcal{F}),\eta\otimes(\kappa\boxtimes
P))$, we can construct a measure structure on $C\times \Omega$ as
follows. Let $\mathcal{E}=\{D\cap (C\times \Omega):
D\in\mathcal{L}\otimes(\mathcal{K}\boxtimes\mathcal{F}) \}$ (which
is a $\sigma$-algebra on $C\times \Omega$), and $\nu$ be the set
function on ${\cal E}$ defined by $\nu(D\cap (C\times
\Omega))=\eta\otimes(\kappa\boxtimes P)(D)$ for any measurable set
$D$ in $\mathcal{L}\otimes(\mathcal{K}\boxtimes\mathcal{F})$. Then,
$\nu$ is a well-defined probability measure on $(C\times \Omega,
{\cal E})$ since the $\eta\otimes(\kappa\boxtimes P)$-outer measure
of $C \times \Omega$ is one. It is obvious that the process $F$ is
$\mathcal{E}$-measurable.

Next, we show that the probability space $(C\times \Omega, {\cal E},
\nu)$ extends the usual product probability space $(C\times
\Omega,\mathcal{T}\otimes \mathcal{F}, \gamma \otimes P)$. Fix any
$Y \in \mathcal{T}$ and $A \in \mathcal{F}$. Then, there is a
measurable set $U \in \mathcal{L}\otimes\mathcal{K}$ such that $Y =
U \cap C$. The rectangle $Y \times A$ is $(U \times A) \cap (C
\times \Omega)$; and hence it belongs to ${\cal E}$. By the
definitions of $\gamma$ and $\nu$, we know that
\begin{eqnarray*}
(\gamma \otimes P) (Y \times A) &=& \gamma (Y) \cdot P(A) = (\eta
\otimes \kappa) (U) \cdot P(A) = (\eta\otimes \kappa \otimes P) (U
\times A) \\
&=& \eta\otimes(\kappa\boxtimes P) (U \times A)= \nu [(U \times A)
\cap (C \times \Omega)] = \nu(Y \times A).
\end{eqnarray*}
Since ${\cal E}$ contains all the rectangles $Y \times A$ for $Y \in
\mathcal{T}$ and $A \in \mathcal{F}$ (which generate
$\mathcal{T}\otimes \mathcal{F}$), it contains $\mathcal{T}\otimes
\mathcal{F}$. Since the probability measures $(\gamma \otimes P)$
and $\nu$ agree on all the rectangles $Y \times A$ for $Y \in
\mathcal{T}$ and $A \in \mathcal{F}$ (which generate
$\mathcal{T}\otimes \mathcal{F}$ and form a $\pi$-system), the
theorem on the uniqueness of measure (p.~404 of \cite{Durrett})
implies that $(\gamma \otimes P)$ and $\nu$ must agree on
$\mathcal{T}\otimes \mathcal{F}$. Therefore, $(C\times \Omega, {\cal
E}, \nu)$ is an extension of $(C\times \Omega,\mathcal{T}\otimes
\mathcal{F}, \gamma \otimes P)$.

In the following, we shall show that  $(C\times \Omega, {\cal E},
\nu)$ is a Fubini extension. Fix any measurable set $E\in
\mathcal{E}$. Then, $E=D\cap (C\times \Omega)$ for some $D\in
\mathcal{L}\otimes(\mathcal{K}\boxtimes \mathcal{F})$.

For each $\omega\in\Omega$, $E_\omega=D_\omega\cap C$. By the Fubini
property associated with $\eta\otimes(\kappa\boxtimes P)$, we have
for $P$-almost all $\omega\in\Omega$, $D_\omega$ is in
$\mathcal{L}\otimes\mathcal{K}$, which means that $E_\omega \in
{\cal T}$ with $\gamma(E_\omega)=\eta\otimes \kappa(D_\omega)$. By
the same Fubini property again, $\eta\otimes \kappa(D_\omega)$ is
$P$-integrable with integral $\int_\Omega \eta\otimes
\kappa(D_\omega)\ud P=\eta\otimes(\kappa\boxtimes P)(D)$. Hence,
$\int_\Omega \gamma(E_\omega)\,\ud P= \eta\otimes(\kappa\boxtimes
P)(D)$, which implies that $\int_\Omega \gamma(E_\omega)\,\ud P =
\nu(E)$ by the definition of $\nu$.

Next we shall prove the other part of the Fubini property associated
with $\nu$ for the measurable set $E\in {\cal E}$. Recall that
$\nu(E) = \eta\otimes(\kappa\boxtimes P)(D)$. By the Fubini property
associated with $\eta \otimes ( \kappa \boxtimes P)$, the function
$P (D_{(l, k)}) $ on $L \times K$ is integrable over $(L\times K,
{\cal L \otimes K}, \eta \otimes \kappa)$ with integral $\eta
\otimes ( \kappa \boxtimes P) (D) = \int_{L \times K}\ P( D_{(l,
k)})\, \ud \,\eta \otimes \kappa$. Since $(L \times K, {\cal U},
{\gamma })$ is an extension of $(L\times K, {\cal L \otimes K}, \eta
\otimes \kappa)$, $P (D_{(l, k)})$ is integrable over  $(L \times K,
{\cal U}, {\gamma })$ with $\int_{L\times K}\ P(D_{(l,k)})\ \ud
\,{\gamma} = \int_{L \times K}\ P( D_{(l, k)})\, \ud \,\eta \otimes
\kappa$. Since $C \in {\cal U}$ with measure $ \gamma (C) =1$ and
${\cal T}$ is the restriction of ${\cal U}$ to $C$, the restriction
of $ P( D_{(l, k)})$ to $C$ is integrable over $(C, \mathcal{T},
\gamma)$ with $\int_C P(D_{(l, k)})\, \ud \gamma = \int_{L\times K}\
P(D_{(l,k)})\ \ud \,{\gamma}$. Since $E_{(l, k)}= D_{(l, k)}$ for
any $(l, k)\in C$, we know that $ P( E_{(l, k)})$ is integrable over
$(C, \mathcal{T}, \gamma)$ with $\int_C P(E_{(l, k)})\, \ud \gamma =
\int_C P(D_{(l, k)})\, \ud \gamma$. By combining all these
equalities together, we obtain that $\nu(E) = \int_C P(E_{(l, k)})\,
\ud \gamma$.

Therefore the indicator function $1_E$ satisfies the Fubini property
for any measurable set $E \in {\cal E}$. The rest of the proof of the Fubini property is the same as in Step 1.
Thus the probability space $(C \times \Omega, {\cal E}, \nu )$ is a
Fubini extension of the usual product probability space $(C \times
\Omega, \mathcal{T} \otimes \mathcal{F}, \gamma \otimes P)$.

\medskip
\nt $\textbf{Step 3.}$ Now let $f: I \times \Omega \rightarrow
[0, 1]$ be another process defined by $f(i,\omega) =  F(\psi^{-1}
(i),\omega)$ for any $(i,\omega)\in I\times \Omega$, where $F$ is
the process in Step 2 and $\psi$ the isomorphism between the
probability spaces $(C, \mathcal{T}, \gamma)$ and $(I, \mathcal{I},
\lambda)$. It is clear that the process $f$ is essentially pairwise
independent and the random variable $f_i$ has uniform distribution
$\mu$ on $[0,1]$ for any $i \in I$.

Given the Fubini extension $(C \times \Omega, \mathcal{E}, \nu)$ of
the usual product probability space $(C \times \Omega,\mathcal{T}
\otimes \mathcal{F}, \gamma \otimes P)$ in Step 2, we can use the
bijection $(\psi, \Id_\Omega)$ from $C \times \Omega$ to $I \times
\Omega$ to construct a $\sigma$-algebra ${\cal W}=\{H \subseteq I \times \Omega: (\psi,
\Id_\Omega)^{-1}(H)\in \mathcal{E}\}$ on $I\times \Omega$, where
$\Id_\Omega$ is the identity map on $\Omega$. Define a probability
measure $\rho$ on $\cal W$ by letting $\rho(H) = \nu [(\psi,
\Id_\Omega)^{-1}(H)]$ for any $H \in {\cal W}$. Therefore, $(\psi,
\Id_\Omega)$ is also an isomorphism between the two probability
spaces $(C \times \Omega, \mathcal{E}, \nu)$ and $(I \times \Omega,
{\cal W},\rho)$. The process $f$ is obviously
$\mathcal{W}$-measurable.

For $S \in {\cal I}$ and $Y \in {\cal F}$, the definition of ${\cal I}$ implies that $\psi^{-1} (S) \in {\cal T}$, and hence $(\psi,
\Id_\Omega)^{-1}(S \times Y)\in \mathcal{T}
\otimes \mathcal{F} \subseteq \mathcal{E}$. Therefore, the definition of ${\cal W}$ implies that $S \times Y \in {\cal W}$. By the definition of $\rho$,
\begin{eqnarray*}
\rho (S \times Y ) &=& \nu [(\psi, \Id_\Omega)^{-1}(S \times Y)] =
\nu [ \psi^{-1} (S) \times Y]
= \gamma \otimes P [ \psi^{-1} (S) \times Y]\\
&=& \gamma [\psi ^ {-1} (S)]\cdot P (Y) = \lambda (S) P (Y) = \lambda \otimes P (S \times Y).
\end{eqnarray*}
The probability measures $\rho$ and $\lambda \otimes P$ agree on all the rectangles $S \times Y$ for $S \in {\cal I}$ and $Y \in {\cal F}$, which generate
$\mathcal{I}\otimes \mathcal{F}$ and form a $\pi$-system. As in Step 2, the
theorem on the uniqueness of measure (p.~404 of \cite{Durrett})
implies that $(\lambda \otimes P)$ and $\rho$ must agree on
$\mathcal{I}\otimes \mathcal{F}$. Therefore, $(I\times \Omega, {\cal
W}, \rho)$ is an extension of $(I\times \Omega,\mathcal{I}\otimes
\mathcal{F}, \lambda \otimes P)$.

Next, we prove the Fubini property associated with $\rho$. As in Step
2, we only prove this property for any measurable set $H \in {\cal
W}$. Let $E = (\psi, \Id_\Omega)^{-1}(H)$; then $E \in {\cal E}$ and  $\nu (E) = \nu [(\psi, \Id_\Omega)^{-1}(H)] = \rho (H) $.

It is obvious that
for any $\omega \in \Omega$, $\psi ^{-1} (H_\omega) = E_\omega$, and $H_\omega = \psi(E_\omega)$. By the definition of $\lambda$, $\lambda
(H _ \omega) = \gamma [\psi ^{-1} (H _ \omega)] = \gamma (E_\omega) $.
By the Fubini property associated with $\nu$ for $E$,
for $P$-almost all $\omega \in \Omega$, $E_ \omega \in {\cal T}$, and thus it follows from the definition of ${\cal I}$ that $H_\omega = \psi(E_\omega) \in {\cal I}$.
By the Fubini property associated with $\nu$ for $E$ again, the $P$-integrable
function $\gamma (E_\omega)$ has integral $\nu (E) = \int_\Omega
\gamma (E_\omega)\, \ud P$. Since $\rho (H) = \nu (E)$, we have
$\rho (H) = \int_\Omega \lambda (H _ \omega)\, \ud P$.

For the other part of the Fubini property associated with $\rho$ for $H$, note that $H
_ i = E _ {\psi ^{-1} (i)}$ for each $i \in I$. By the Fubini property of $\nu$ for $E$, there is a set $T \in \mathcal{T}$ with $\gamma(T) =1$ such that for any $(l ,k) \in T$, $E _ {(l, k)} \in \mathcal{F}$. Hence, $\psi(T) \in \mathcal{I}$ with $\lambda(\psi(T)) =1$, and for each $i \in \psi(T)$, $H_i = E _ {\psi ^{-1} (i)} \in \mathcal{F}$. By the Fubini property of $\nu$ for $E$ and the formula for changing variables,
\begin{eqnarray*}
 \nu (E) &=& \int _ C P (E _ {(l, k)}) \, \ud \gamma = \int_T P (E _ {(l, k)}) \, \ud \gamma = \int_{\psi(T)} P(E _ {\psi ^{-1} (i)}) \, \ud \lambda \\
&=& \int_{\psi(T)} P(H_i) \, \ud \lambda = \int_I P(H_i) \, \ud \lambda.
\end{eqnarray*}
Since $\rho (H) = \nu (E)$, we have
$\rho (H) = \int_I P(H_i) \, \ud \lambda$.

Therefore, $(I \times \Omega, {\cal W},\rho)$ is a Fubini extension. As in Definition \ref{D:Fubini extension}, we denote $(I \times
\Omega, {\cal W},\rho)$ by $(I \times \Omega, \mathcal{I} \boxtimes
\mathcal{F}, \lambda \boxtimes P)$.  $\Box$

\ms

In the following lemma, we restate Proposition 5.3 of~\cite{Sun-06}.

\begin{lemma}\label{L:vary distributions.}
Let $(I\times\Omega,\mathcal{I}\boxtimes\mathcal{F},\lambda\boxtimes P)$ be a Fubini extension, and $X$ a complete separable metric space.
Assume that there exists
an essentially pairwise independent process $f$ from $I\times
\Omega$ to $[0,1]$ such that $f$ is
$\mathcal{I}\boxtimes\mathcal{F}$-measurable, and for each $i \in
I$, the distribution of the random variable $f_i$ is the uniform
distribution $\mu$ on $[0, 1]$.\fn{Such a Fubini extension is called a \emph{rich
product probability space} in \cite{Sun-06}.} Then,
for any measurable mapping $\varphi$ from $(I, {\cal I},
\lambda)$ to the space ${\cal M}(X)$ of Borel probability
measures on $X$, there is a ${\cal I} \boxtimes
{\cal F}$-measurable process $g$ from $I \times \Omega$ to $X$ such
that the random variables $g_i$ are essentially pairwise
independent, and the distribution $P g_i^{-1}$ is the given
distribution $\varphi(i)$ for $\lambda$-almost all $i \in I$.
\end{lemma}

 \nt \textbf{Proof of the
Theorem~\ref{T: main theorem}:} It is now obvious that Theorem~\ref{T: main theorem} follows from Proposition \ref{P:rich product space.} and Lemma \ref{L:vary distributions.}. $\Box$

\end{document}